\documentclass[12pt,reqno]{amsart}

\usepackage{color}
\usepackage{amsmath, amsthm, amssymb}
\usepackage{amsfonts}
\usepackage[ansinew]{inputenc}
\usepackage[dvips]{epsfig}
\usepackage{graphicx}
\usepackage[english]{babel}
\usepackage{hyperref}

\usepackage{cite}
\usepackage{graphicx}
\usepackage{amscd}
\usepackage{color}
\usepackage{bm}
\usepackage{enumerate}

\usepackage{verbatim}
\usepackage{hyperref}
\usepackage{amstext}
\usepackage{latexsym}
\theoremstyle{plain}
\newtheorem{Theorem}{Theorem}[section]

\newtheorem{Corollary}[Theorem]{Corollary}

\newtheorem{Remark}[Theorem]{Remark}

\numberwithin{Theorem}{section} \numberwithin{equation}{section}

\def\square{\vbox{
\hrule height .4pt \hbox{\vrule width .4pt height 7pt \kern 7pt
\vrule width .4pt} \hrule height .4pt }}
\def\QED{\hfill {$\square$}\goodbreak \medskip}

\newcommand{\average}{{\mathchoice {\kern1ex\vcenter{\hrule height.4pt
width 6pt depth0pt} \kern-9.7pt} {\kern1ex\vcenter{\hrule height.4pt
width 4.3pt depth0pt} \kern-7pt} {} {} }}

\def\R{\mathbb{R}}

\def\div{\text{div}}

\renewcommand{\a }{\alpha }

\newcommand{\e }{\varepsilon }

\newcommand{\n }{\nabla }
\newcommand{\vp }{\varphi }

\newcommand{\s }{\sigma }

\newcommand{\z }{\zeta}

\renewcommand{\o }{\omega }
\renewcommand{\O }{\Omega }

\newcommand{\ov}{\overline}

\newcommand{\be}{\begin{equation}}
\newcommand{\ee}{\end{equation}}

\newcommand{\de}{\partial}

\newcommand{\ti}{\widetilde}

\newcommand{\calH }{\mathcal{H}}

\newcommand{\calC }{\mathcal{C}}

\newcommand{\N}{\mathbb{N}}

\newcommand{\Z}{\mathbb{Z}}

\newcommand{\cG}{{\mathcal G}}

\renewcommand{\epsilon}{\varepsilon}

\begin{document}

\title[]{Unbounded periodic constant mean
curvature graphs on  Calibrable Cheeger Serrin domains}


\author{Ignace Aristide Minlend}
\address{Department of Quantitative Techniques, Faculty of Economics and Applied Management, University of Douala, BP. 2701 - Douala.}
\email{ignace.a.minlend@aims-senegal.org}




\begin{abstract}
We prove a general result characterizing  a specific class of Serrin
domains as supports of unbounded and periodic constant mean
curvature graphs. We apply this result to prove the existence of a
family of unbounded periodic constant mean curvature graphs, each
supported by a Serrin domain and intersecting its boundary
orthogonally, up to a translation. We also show that the underlying
Serrin domains are calibrable and Cheeger in a suitable sens, and
they solve the $1$-Laplacian equation.
\end{abstract}

\maketitle \textbf{Keywords}: Overdetermined problems, Cheeger sets,
calibrable sets, Serrin domains, mean curvature. \\


\section{Introduction and main result}

In the more recent literature, a domain $\Omega$ in the Euclidean
space  $\mathbb{R}^{N}$, $N\geq 2$ is called {\em Serrin domain} if
it admits a {\em positive} solution to the following so called
Serrin's overdetermined problem
\begin{align}
 &-\Delta u=1 \qquad \text{ in $\Omega$}\label{eq:Pro-1}\\
&u=0, \quad \partial_\nu u= -\beta\quad \text{on $\partial
\Omega$,}\label{eq:Pro-2}
\end{align}
where  $\nu$ is the unit outer normal on  $\partial \O$ and $\beta$
is a positive constant. In \cite{Serrin} Serrin proved that the only
$C^2$ bounded domains $\O$ where problem
(\ref{eq:Pro-1})-(\ref{eq:Pro-2}) admits a solution are balls. This
parallels the result of Alexandrov \cite{Alexandrov} that smooth
bounded domains with constant boundary mean curvature are balls.

In general, it  is a challenging question  to relate  constant mean
curvature (CMC) surfaces  with domains where Serrin's problem
(\ref{eq:Pro-1})-(\ref{eq:Pro-2}) is solvable. In lower dimension
however, Del Pino, Pacard and Wei \cite{del-Pino-pacard-wei}
succeeded to provide  examples for domains whose boundary is close
to large dilations of a given CMC surface where Serrin’s
overdetermined problem is solvable. In \cite{FarinaVal}, A. Farina
and  E. Valdinoci, and  independently H. Berestycki, L. A.
Caffarelli, and L. Nirenberg \cite{BCNI} also studied the conditions
under which an epigraph domain solving an overdetermined problem
must be a half-space. Observe in this case that the boundary of $\O$
is then a minimal surface.

Under further assumptions, one might also expect a Serrin domain to
support a graph with CMC. For instance if a Serrin domain $\O$ is
bounded, it is known from our result in \cite{I.A.M} that $\O$ is
uniquely self-Cheeger set. Thanks to E. Guisti \cite{Giusti}, there
exists a unique (up to a constant) bounded solution $w\in C^{2}(\O)$
to the prescribed mean curvature equation
\begin{align}\label{eq:mecue}
\begin{cases}
-\div\dfrac{\n w}{\sqrt{1+|\n w|^2}} =h(\O)& \qquad \textrm{ in } \O \vspace{3mm}\\
-\dfrac{\n w}{\sqrt{1+|\n w|^2}}\cdot \eta=1  & \qquad \textrm{ on }
\de\O,
\end{cases}
\end{align}
where $\eta$ is the unit  outer normal to $\de\O$ and $h(\O)$ is the
Cheeger constant of $\O$.  That is $\O$ supports a graph with
constant mean curvature equals to $h(\O)$, and this graph intersects
perpendicularly the horizontal plane, up to a translation.

In this paper, we examine the analogue of this result on a specific
class of unbounded  Serrin domains. Namely  we aim to prove  the
existence of constant mean curvature graphs $w$ satisfying
\eqref{eq:mecue}, with the property that they are supported by
unbounded and periodic calibrable Cheeger Serrin domains $\O$. The
class of  Serrin domains under consideration are domains
$\O\subset\R^N= \R^n \times \R^m$ that are radial  and bounded in
the variable $z\in \mathbb{R}^{n}$, periodic and symmetric in the
variables $t_1,\dots t_m$ in $\mathbb{R}^m$. A prototype of these
domains was already constructed by M. M. Fall, T. Weth and the
author in  \cite[Theorem 1.1]{Fall-MinlendI-Weth}.

The notion of Cheeger sets we are dealing with in this work
generalizes the classical one in \cite{Leon} and involves the
definition of relative perimeter. We let $S$ and $\Omega$ be open
subsets of $\R^N$. For a given Borel set $A$ in $\R^N$, the relative
perimeter of $A$ in $S$ is given by
$$
P(A,S):=\sup\left\{\int_{A} \div(\xi)\, dx\,:\, \xi\in
C^\infty_c(S;\R^{N}),\, |\xi|\leq 1\right\}.
$$ The perimeter, $P(A)$, of a Borel set $A \subset \R^N$
is the relative perimeter   $P(A, \R^N)$ of $A$ in $\R^N$. If a set
$A$ is of finite perimeter, denoting by $\partial^{*}A$ the reduced
boundary of $A$, then by De Giorgi's structure Theorem \cite[Theorem
2.2]{A.FigalliandGe} (see also \cite[Remark on p.161]{Hofmann} or
\cite{AmbrosioPallara}), $P(A, \O) =H^{N-1}(\partial^{*}A\cap \O).$
In addition, we have the Gauss-Green formula
\begin{equation}\label{The div}
\int_{A} \div(\xi)\,\textrm{d}x=\int_{\partial^{*}A} \langle
\xi,\nu_A \rangle \textrm{d} H^{N-1}\quad \textrm{for all}
\quad\xi\in C^1_c(\R^{N};\R^{N}).
\end{equation}
Here and in the following, $H^{N-1}$ denotes the $(N-1)$-dimensional
Hausdorff measure. The Cheeger constant of an open set $\Omega$
relative to $S$ is by definition \be \label{eq:def_Cheeg} h(\O,S):=
\inf_{A\subseteq \O \cap S}\frac{P(A,S)}{|A|}, \ee where the infimum
is taken over Borel subsets $A \subseteq \Omega \cap S$ with  finite
perimeter. If this constant is attained by some Borel subset $A
\subset \O \cap S$ with finite perimeter, then $A$ will be called a
Cheeger set of $\O$ relative to $S$. If  $A=\Omega \cap S$ attains
the constant $h(\O,S)$ in \eqref{eq:def_Cheeg}, we say that $\O$ is
self-Cheeger relative to $S$. Moreover, if any Cheeger set in  $\O
\cap S$ is equal to $\O\cap S$ up to a  set with zero  Lebesgue
measure, we say that $\O$ is uniquely self-Cheeger relative to $S$.

Related to the notion Cheeger sets is the concept calibrable sets
defined in \cite{Alter}, where convex bounded calibrable sets were
characterized. In the same fashion, we introduce the notion of
relative calibrable sets. We say that a set $E$ is calibrable
relative to S, if there  exists a vector field
 $\xi\in L^{\infty}(E, \R^N)$,  $|\xi|_{\infty} \leq 1$ such that
\begin{align}\label{eq:Calibrable-unique}
 -\div \xi=\lambda_E 1_{E} \quad \textrm{in} \quad
 \mathcal{D}'(\R^N),\quad \xi\cdot \nu=-1 \quad \textrm{on} \quad \de^{*}E\cap
 S,\quad
\xi\cdot \nu=0 \quad \textrm{on} \quad  E\cap \de S,
\end{align}
where  $1_{E}$ denotes the characteristic function of $E$. Observe
that if $E$ is calibrable with respect to $S$, then
$\lambda_E=\frac{P(E,S)}{|E|}$. With these notations and
definitions, we can now state our first result.

\begin{Theorem}\label{eq:CMHVdCC}
Let $\O\subset\R^N= \R^n \times \R^m$ be a radial and bounded domain
in the variable $z\in \mathbb{R}^{n}$, $2\lambda$-periodic and
symmetric in the directions  $t_1,\dots t_m$ in $\mathbb{R}^m$.
Assume that the overdetermined problem
(\ref{eq:Pro-1})-(\ref{eq:Pro-2}) admits a solution $u \in
C^{2,\a}(\ov{\O})$.

Then there a function $w\in C^\infty(\O)$ which is radial in  $z\in
\mathbb{R}^{n}$, even and $2\lambda$-periodic in the variables
$t_1,\dots t_m$ and such that
\begin{align}\label{eq:w_s-ok1}
\begin{cases}
\displaystyle -\div\frac{\n w }{\sqrt{1+|\n w |^2}}=\frac{1}{\beta} & \quad \textrm{in} \quad \O  \vspace{3mm}\\
 \displaystyle - \frac{\n w}{\sqrt{1+|\n w|}}\cdot \eta=1  & \quad \textrm{on} \quad\de \O
\end{cases}
\end{align}
Furthermore, the function $w$ is unique (up to additive constant) in
the class of even and $2\lambda$-periodic functions in the variables
$t_1,\dots, t_m$.
\end{Theorem}
 \begin{Remark}
We do not know that $w\in C(\ov{\O})$ and that $w|_{\de \O}\equiv
Const$. This would implies that the graph of $w$ is a half of
Delaunay hypersurface \cite{Delaunay}.
\end{Remark}
Notice that  the existence of solution to \eqref{eq:mecue}  as
stated in \cite{Giusti} is strongly related to the fact that $\O$ is
uniquely self-Cheeger set. It is known from our result in
\cite{I.A.M} that any bounded Serrin domain is uniquely self-Cheeger
set. For  Serrin domains consider in Theorem \ref{eq:CMHVdCC}, we
show that they are \emph{uniquely self-Cheeger} and calibrable
relatively to a tubular domain in $\mathbb{R}^N$. In addition, they
are support of a solution to the $1$-Laplacien equation.

\begin{Theorem}\label{cor:Cheeger}
Let $\O\subset\R^N= \R^n \times \R^m$ be a radial and bounded domain
in the variable $z\in \mathbb{R}^{n}$, $2\lambda$-periodic and
symmetric in the directions  $t_1,\dots t_m$ in $\mathbb{R}^m$.
Assume that the overdetermined problem
(\ref{eq:Pro-1})-(\ref{eq:Pro-2}) admits a solution $u \in
C^{2,\a}(\ov{\O})$ and let $a, b \in \lambda \Z^m$ with $a_i < b_i$
for $i=1,\dots,m$.

Then $\O$  is uniquely self-Cheeger relative to the set
\begin{equation}
\label{eq:defsalphabeta}
 S_a^b:= \R^n \times (a_1, b_1) \times \dots \times (a_m, b_m)
\end{equation}
with corresponding relative Cheeger constant $h(\Omega ,S_{a}^b)=
\dfrac{1}{\beta}.$

Moreover, \be \label{eq:1-Lapl} -\textrm{div}\left(\frac{\n 1_{\O}
}{|\n 1_{\O}|} \right)=\frac{1}{\beta} 1_{\O} \ee and  $\O$ is
calibrable relative to $S_a^b$.
\end{Theorem}

\begin{Remark}\label{Remark}
As a remark, the first statement in Theorem \ref{cor:Cheeger} is the
result of \cite[Corollary 1.2]{Fall-MinlendI-Weth}, where the
Cheeger constant was defined taking the infimum over subsets $A
\subset \Omega \cap S$ with Lipschitz boundary. The same conclusion
in \cite[Corollary 1.2]{Fall-MinlendI-Weth} holds even for current
definition of the Cheeger constant in \eqref{eq:def_Cheeg}. Indeed
working with \eqref{eq:def_Cheeg}, the main challenge in the proof
of the first part of Theorem \ref{cor:Cheeger} is to show the {\em
uniqueness} of $\O$ as a self-Cheeger relative $S_a^b$. But scanning
\cite[Section 5]{Fall-MinlendI-Weth}, this follows from the last
part of the proof Theorem 1.2 in \cite{I.A.M} after replacing $\O$
by $\O\cap S_a^b$.
\end{Remark}

\begin{Corollary}\label{eq:CMHVdC}
There exists a family of unbounded constant mean curvature graphs
$(\Gamma_s)_{s\in(-\e_0, \e_0)}$, each supported by a calibrable
Cheeger Serrin domain $E_s$ and intersecting its boundary
orthogonally (up to a translation).

For every $s\in (-\e_0, \e_0)$, the graph $\Gamma_s$ is
$2\pi$-periodic in some variables and radial in the others, and has
mean curvature equal to $\frac{n}{\lambda_s}$, with  $\lambda_s>0.$
\end{Corollary}

The existence of Serrin domains stated in Corollary \ref{eq:CMHVdC}
relies on our recent result in \cite[Theorem
1.1]{Fall-MinlendI-Weth}, where we constructed a family
$(E_s)_{s\in(-\e_0, \e_0)}$ of unbounded periodic Serrin domains of
the form
$$
E_s:=\left\{ (z,t)\in \mathbb{R}^{n}\times\mathbb{R}^m\,:\,
|z|<\phi_s(t)  \right\}\subset\R^N,
$$
where $N=n+m$ and $\phi_s: \R^m \to (0,\infty)$ is an even and $2\pi
\Z^m$-periodic function. We  note that, each of $E_s$ solves the
overdetermined (\ref{eq:Pro-1})-(\ref{eq:Pro-2}) with
$\beta=\dfrac{\lambda_s}{n}$.\\

From Theorem \ref{cor:Cheeger}, Serrin domains $E_s$ are calibrable
and Cheeger. The proof of Corollary \ref{eq:CMHVdC} is then complete
once we prove the more general result in Theorem \ref{eq:CMHVdCC},
namely the existence of a solution to \eqref{eq:w_s-ok1}.

The proof of Theorem \ref{eq:CMHVdCC} is inspired by \cite{Giusti-1}
and \cite{Giusti}. We first apply [\cite{Fall-MinlendI-Weth}, Lemma
5.1] to the overdetermined problem (\ref{eq:Pro-1})-(\ref{eq:Pro-2})
and show the conditions (1.4) and (1.4) of [\cite{Giusti-1}, Section
1.C] are fulfill by the minimization problem associated with
\eqref{eq:mecue}. The result of [\cite{Giusti-1}, Theorem 2.1] then
yields the existence of a function $v_{\e}\in BV(F_{\e})\cap
C^{2}(F_{\e})$ solution to \eqref{eq:Solvk}. By compactness, the
sequence $(v_{\frac{1}{n}})_{n\in\N}$ converges to a generalized
solution $v: F \to [0,+\infty]$, which  we prove to be a classical
solution by showing it is bounded. The function $v$ is even in the
variables $t_1,\dots,t_m$  and this allows to derive the second
equation of \eqref{eq:Solvko-ok}. To finish the proof, we make the
shift $\ti{v}(z,t):=v(z,t_1-2\lambda, \dots,t_m-2\lambda )$, and
finally deduce in \eqref{eq:dew} the existence of a function $w\in
C^\infty(\O)$, radial in the $z$ variable and $2\lambda$-periodic in
the variables $t_1,\dots, t_m$, which in addition satisfies
\eqref{eq:w_s-ok1}.\\

\bigskip
\noindent \textbf{Acknowledgement:} The author is grateful to  Pr.
Mouhamed Moustapha Fall and Pr. Tobias Weth for their helpful
suggestions and insights throughout the writing of this paper. Part
of this work was achieved when the author was awarded the inaugural
"Abbas Bahri Excellence Fellowship" at Rutgers University of
New-Jesey, USA. He wishes to thank  the Department of Mathematics
for the hospitality.

\section{Proof of Theorem \ref{cor:Cheeger}}\label{ProofCo}
From Remark \ref{Remark}, we only need to prove the last statement
of Theorem \ref{cor:Cheeger}.

We start by explaining the meaning of \eqref{eq:1-Lapl}(see also
\cite{SchevenSchmidt}).

Let $\O$ be an open set of $\R^{N}$ and $u \in L^1(\O)$, the total
variation of $u$ in $\O$ is defined by
$$
|\n u|(\O)=\sup\left\{\int_{\O} u\div(\xi)\, dx\,:\, \xi\in
C^\infty_c(\O,\R^{N}),\, |\xi|\leq
 1\right\}.
$$
The space of functions with bounded variation in an open $\O$ is
denote by  $BV(\O)$. For a given Borel set $A$, it is plain that
$P(A,S)= |\n 1_{A}|(S)$.

Let   $f\in L^1_{loc}(\O)$. A  function $v\in BV_{loc}(\O)\cap
L^\infty_{loc}(\O)$ solves the equation
$$
-\div \frac{\n v}{|\n v |}= f \qquad \textrm{ in }\O
$$
if there exists a sub-unit vector field $\xi \in W^{1,1}_{loc}( \O,
\R^{N})$ such that \be\label{eq:pair-nu} (\xi,\n v)=|\n v|\quad
\textrm{ as Radon measures  in $\O$ } \ee and \be\label{eq:pair-nu1}
-\div \xi =f  \quad  \textrm{ in the sense of distribution in $\O$}.
\ee We recall the Radon measure $(\xi, \n v)$ in \eqref{eq:pair-nu}
is defined (see Anzelotti \cite{Anz}) thanks to Riesz representation
theorem via distributional sense  as \be (\xi, \n v)(\vp)=-\int_{\O}
v \vp  \div \xi\, dx -\int_{\O} v \xi\cdot \n\vp\, dx \qquad
\textrm{ for every $\vp\in C^\infty_c(\O)$}. \ee We remark,   that
if $\O$ is a Lipshcitz open set and $v=1_{\O}$, we have
\begin{equation}\label{eqradmeas}
(\xi, \n v) =-\xi\cdot \nu_\O \calH^{N-1}|_{\de \O}.
\end{equation}
by integration by parts. In addition, it is well know  from the
divergence theorem that
\begin{equation}\label{eqhausmeas}
 |\n 1_{\O}|=  \calH^{N-1}|_{\de \O},
\end{equation}
where $\calH^{N-1}$ is the $(N-1)$-dimensional Hausdorff measure. \\

Since $u=0$ on $\de \O$, we have $\n u=-|\n u|\eta$ on $\de\O$,
where $\eta$ denotes the outer unit normal vector field of $\de \O$.
We put
$$
\xi = \frac{1}{\beta} \n  u\quad \textrm{and } v= 1_{\O}.
$$
It is plain that the solution $u$ of problem
(\ref{eq:Pro-1})-(\ref{eq:Pro-2}) satisfies the gradient estimate
\begin{equation}\label{eqestigraaa}
|\n u|< \beta  \quad \textrm{ in }\O.
\end{equation}
This follows from [\cite{Fall-MinlendI-Weth}, Lemma 5.1], see also
[\cite{I.A.M}, Section 2] for the corresponding estimate in compact
Riemannian manifolds.

Now using \eqref{eqradmeas}, \eqref{eqhausmeas} and
\eqref{eqestigraaa}, we find that  \eqref{eq:pair-nu} and
\eqref{eq:pair-nu1} are satisfied with
$$f=\frac{1}{\beta}~1_{\O}$$
and this yields \eqref{eq:1-Lapl}.

 We now  end the proof by showing that the the set
$\O$ is calibrable relative to $S_a^b$.

We observe that we have just constructed the vector field
$\xi=\frac{1}{\beta} \,\n u$ which by \label{eqestigra} clearly
satisfies $|\xi|\leq 1$ by \eqref{eqestigraaa}.

From the first part of the proof, we also have
\begin{align}\label{eq:Cable-unique}
-\div \xi=\frac{1}{\beta}\, 1_{\O}=h(\O,S_{a}^b)1_{\O} \\
\xi=-\eta \quad \textrm{on} \quad \de \O\cap S_a^b,\nonumber
\end{align}
where for $a, b \in \lambda\Z^m$ with $a_i < b_i$ ($i=1,\dots,m$),
the set $S_a^b$ is defined as in (\ref{eq:defsalphabeta}). We
emphasize that the boundary $\partial S_a^b$ can be decomposed into
a disjoint union $\partial S_\tau= K \cup S^1 \cup \dots \cup S^m$,
where
$$
S^i:= \R^n \times \{t \in \R^m  \;:\; t_i \in \{a_ i, b_i\},\; t_j
\in (a_ j , b_j) \; \text{for $j \not = i$}\}\qquad \text{for
$i=1,\dots,m$,}
$$
and $K$ has zero $(N-1)$-dimensional Hausdorff measure. Since by
hypothesis $\O$ is $2\lambda$-periodic and symmetric in
$t_1,\dots,t_m$, the solution $u$ of problem
(\ref{eq:Pro-1})-(\ref{eq:Pro-2}) is $2\lambda$-periodic and even in
$t_1,\dots,t_m$ which leads to
$$
\frac{\partial u}{\partial {t_i}}  \equiv 0 \qquad \text{on $S^i$
for $i=1,\dots,m$}.$$ Observe that the outer unit normal to $S^i$
coincides with $(0,e_i)$ or $(0,-e_i)$, where $0 \in \R^n$ and $e_i$
denotes the $i$-th coordinate vector in $\R^m$ and we have
\begin{align}\label{pa-der}
  \xi\cdot e_{i}=\frac{1}{\beta}\frac{\partial u}{\partial {t_i}}
\equiv 0 \qquad \text{on $S^i$ for $i=1,\dots,m$.}
\end{align}
Hence, following \cite{Alter} and our definition
\eqref{eq:Calibrable-unique}, we can say from
\eqref{eq:Cable-unique} and \eqref{pa-der} that the set $\O$ is
\textit{calibrable} relative to $S_a^b$. \QED

\section{Proof of Theorem \ref{eq:CMHVdCC}}
We prove that under the hypotheses on the Serrin domain $\O$ in
Theorem \ref{eq:CMHVdCC}, there a function $w\in C^\infty(\O)$ which
is even  and $2\lambda$-periodic in the variables $t_1,\dots t_m$,
radial in $z\in \R^n$, and such that \eqref{eq:w_s-ok1} holds.\\

Following [\cite{Giusti-1}, Section 1], we consider the set
$$A_{1}:=(\O \cup \{(z,t)\in \R^n\times \R^m: \textrm{dist}(z,
\partial \O)<1\})\cap S, \quad  \textrm{where}  \quad S:=\R^n \times (-2\lambda, 2\lambda)^m.$$  For
$\e>0$, we also define  the Libschitz domains
$$\O_{\e}:=\{(z,t)\in \O: \quad \textrm{dist}(z, \partial \O)>\e
\}, \quad F_{\e}:=\O_{\e}\cap S,\quad F:= \O \cap S. $$  We  further
fix  notations for boundaries:
$$
\de_1 F_{\e}=\de F_{\e}\cap A_{1}\quad\textrm{and}\quad \de_2
F_{\e}= \ov{F_{\e}}\cap  \de A_{1}.
$$
Let  $B$ be a Borel set in  $F_{\e}$. We use \eqref{The div}
\eqref{pa-der}  and integrate (\ref{eq:Pro-1})-(\ref{eq:Pro-2}) over
$B$ to get
$$
 |B|\leq\sup_{\overline{F_{\e}}}|\n u|\, |\de^* B \cap S| .
$$
We  put $$c_{\e}:=1-\frac{1}{\beta}\sup_{\ov{F_{\e}}}|\n u|.$$ For
every  $\e>0$, $\ov{F_{\e}}\subset \O $ and by \eqref{eqestigraaa},
$c_{\e}>0$ and we  have for any Borel set $B\subset F_{\e} $, \be
 \frac{1}{\beta} |B|\leq (1-c_{\e})|\n 1_B|(S).
\ee
Thanks to [\cite{Giusti-1}, Theorem 2.1], we can find a function
$v_{\e}\in BV(F_{\e})\cap C^{2}(F_{\e})$ solution to
\begin{align}\label{eq:Solvk}
\begin{cases}
\displaystyle -\div\frac{\n v_{\e}}{\sqrt{1+|\n v_{\e}|^2}}=\frac{1}{\beta} & \qquad \textrm{ in } F_{\e} \vspace{3mm}\\
\displaystyle v_{\e}=0  & \qquad \textrm{ on }  \de_1 {F_{\e}} \vspace{3mm}\\
\displaystyle \frac{\n v_{\e}}{\sqrt{1+|\n v_{\e}|^2}}\cdot e_i=0,
\quad i=1,..., m. & \qquad \textrm{ on } \de_2  F_{\e}.
\end{cases}
\end{align}
The second condition is understood in the sense of trace and the
third condition in the following sense \be\label{eq:Neum-k} \lim_{\e
\to 0}\int_{\de_2 F_{\e}}  \frac{\n v_{\e}}{\sqrt{1+|\n
v_{\e}|^2}}\cdot e_i \, d\s_\e=0;  \quad i=1,\dots, m. \ee Hence the
solution is unique up to a additive constant.

Let $\ti{v}_{\e}(z,t):=v_{\e}(z,-t)$ then the function $\ti{v}_{\e}
$ satisfies similar equation as $v_{\e}$ and thus ${v}_{\e}$ is even
in the  variables $t_1,\dots t_m$ and  we have $\de_{t_i}
v_{\e}(z,0)=0, \quad \textrm{for} \quad i=1,\dots, m.$

By a similar argument, we can obtain $v_{\e}(z_1,\dots,-z_i,\dots,
z_n,t)=v_{\e}(z,t)$ and thus
$$
v_{\e}(z,t)=\ov{v^{\e}}(|z|,t).
$$
We recall that   $v_{\e}$ minimizes the functional \be
\cG(w)=\int_{A_{1}}\sqrt{1+|\n w|^2}-\frac{1}{\beta}\int_{ A_{1}}w\,
dx \ee in the class \be \calC:=\{w\in BV(A_{1})\,:\, w=0 \textrm{ in
}{A_{1}}\setminus F_{\e}\}, \ee where
$$
\int_{E}\sqrt{1+|\n w|^2}=\sup\left\{\int_{E} [u\div(\z)+\z_{N+1}]\,
dx\,:\, \z\in C^\infty_c(E;\R^{N+1}),\, |\z|\leq 1\right\}.
$$
Since $\cG(|w|)\leq \cG(w)$, we can assume that $v_{\e}\geq 0$ in
$F_{\e}$.

By compactness (see \cite{Miranda-1977}) the sequence
$(v_{\frac{1}{n}})_{n\in\N}$ converges to a generalized solution $v:
F\to [0,+\infty]$. We now show that $v$ is bounded.

We define  $P:=\{v=+\infty\}$.
 By simple arguments, using also the
co-area formula (see also \cite{Miranda-1964}), one has that $P$
minimizes the functional
 $$
E\mapsto  P(E,S)-\frac{1}{\beta} |E\cap S|.
 $$
In addition, by Theorem \ref{cor:Cheeger}, $
P(\O,S)=\frac{1}{\beta}|\O \cap S|.$
Hence using similar argument as
in [\cite{Giusti}, Lemma 1.2], we deduce that the set $P=\emptyset$.
We then conclude that $v$ is a classical solution and moreover since
it is even in  the variables $t_1,\dots,t_m$ (as pointwise limit of
even functions), we get $\de_{t_i} v(z,0)=0, \quad i=1,\dots, m$.

Now by integrating over $F_{\e}$  we have \be - \int_{\de_1
F_{\e}}\frac{\n v}{\sqrt{1+|\n v|^2}}\cdot \eta^{\e}\, d\s_\e =
\frac{1}{\beta}\frac{|F_{\e}  | }{|\de_1 F_{\e}|}|\de_1 F_{\e}| \ee
so that \be \lim_{\e\to 0} \int_{\de_1 F_{\e}}\left\{-\frac{\n
v}{\sqrt{1+|\n v|^2}}\cdot \eta^\e-1\right\}\, d\s_\e=\lim_{\e\to 0}
|\de_1 F_{\e}|\biggl[\frac{1}{\beta}\frac{|F_{\e}|}{|\de_1 F_{\e}|}
-1\biggl]=0. \ee We then conclude that $v$ solves uniquely (up to an
additive constant) the equation
\begin{align}\label{eq:Solvko-ok}
\begin{cases}
\displaystyle -\div\frac{\n v }{\sqrt{1+|\n v |^2}}=\frac{1}{\beta} & \qquad \textrm{ in } F \vspace{3mm}\\
 \displaystyle - \frac{\n v }{\sqrt{1+|\n v |^2}}\cdot \eta=1  & \qquad \textrm{ on } \de_1  F \vspace{3mm}\\
 \displaystyle \frac{\n v}{\sqrt{1+|\n v|^2}}\cdot e_i=0  & \qquad \textrm{ on } \de_2  F , \quad i=1,\dots, m.
\end{cases}
\end{align}
We now make the shift $\ti{v}(z,t):=v(z,t_1-2\lambda, \dots,
t_m-2\lambda )$ for $(z,t)\in W \cap F$,  where   $W:=\R^{n}\times
(0, 2\lambda)^m$.

By a direct computation, we have that  both $\ti{v}$ and $v$ satisfy
the equation
 \begin{align}\label{eq:Soltivko-ok}
\begin{cases}
\displaystyle -\div\frac{\n \o }{\sqrt{1+|\n \o |^2}}=\frac{1}{\beta} & \qquad \textrm{ in } F\cap W  \vspace{3mm}\\
 \displaystyle - \frac{\n \o}{\sqrt{1+|\n \o|}}\cdot \eta=1  & \qquad \textrm{ on } \de_1  (F\cap W)   \vspace{3mm}\\
 \displaystyle \frac{\n \o}{\sqrt{1+|\n \o|^2}}\cdot e_i=0, \quad i=1,\dots, m & \qquad \textrm{ on } \de_2 ( F\cap W ).
\end{cases}
\end{align}
 By uniqueness, we deduce that $\ti{v}(z,t)=v(z,t_1-2\lambda, \dots, t_m-2\lambda )=v(z,t)+c$,  for some constant $c\in\R$.
 Since  $v(z,\lambda, \dots, \lambda)=v(z,-\lambda, \dots,-\lambda)$, we infer that $c=0$ and thus
\begin{equation}\label{eq:vonW}
v(z, t_1-2\lambda, \dots, t_m-2\lambda )=v(z,t) \quad\textrm{ for
every}\quad (z,t)\in W \cap F.
\end{equation}
This implies that $v\in C^\infty (F\cup \de S)$. For $k\in \Z$,
define the sets $$S_k:=\R^n\times (2k\lambda, 2(k+1)\lambda)^m.$$ If
$(z,t)\in \O\cap S_k,$ then $(z,t_1-2k\lambda, \dots, t_m-2k\lambda
)\in W \cap F$ and from \eqref{eq:vonW} we can define
\begin{equation}\label{eq:dew}
w(z,t):=v(z,t_1-2k\lambda, \dots, t_m-2k\lambda )\quad
\textrm{for}\quad (z,t)\in \O\cap S_k.
\end{equation}
It is then clear  from \eqref{eq:dew} and \eqref{eq:vonW} that the
function $w$ is $2\lambda$-periodic in the variables $t_1,\dots,
t_m$ and satisfies \eqref{eq:w_s-ok1}. In addition  $w\in
C^\infty(\O)$ and is unique, up to a constant, in the class of
functions which are even and $2\lambda$-periodic in the variables
$t_1,\dots, t_m$.\QED

\end{document}